\documentclass[reqno, a4paper]{amsart}

\usepackage[T1]{fontenc}
\usepackage[british]{babel}
\usepackage{graphicx}
\usepackage{amsfonts}
\usepackage{amssymb}
\usepackage{amsthm}
\usepackage{amsmath}
\usepackage[all]{xy}
\usepackage{hyperref}
\usepackage{mathbbol}
\usepackage{mathabx}

\theoremstyle{plain}
\newtheorem{theorem}[subsection]{Theorem}
\newtheorem{lemma}[subsection]{Lemma}
\newtheorem{proposition}[subsection]{Proposition}

\theoremstyle{definition}
\newtheorem{definition}[subsection]{Definition}
\newtheorem{remark}[subsection]{Remark}

\newtheorem{example}[subsection]{Example}

\newcommand{\defn}{\textbf}

\newcommand{\UCE}{condition~{\rm (UCE)}}
\newcommand{\UCEspecial}{condition~(UCE)}
\newcommand{\UCEcapital}{Condition~(UCE)}

\newcommand{\meet}{\ensuremath{\wedge}}
\newcommand{\comp}{\raisebox{0.2mm}{\ensuremath{\scriptstyle{\circ}}}}
\newcommand{\To}{\Rightarrow}
\newcommand{\coker}{\ensuremath{\mathrm{coker\,}}}
\newcommand{\del}{\ensuremath{\partial}}
\newcommand{\Ker}{\ensuremath{\mathrm{Ker}}}
\renewcommand{\H}{\ensuremath{\mathrm{H}}}
\newcommand{\U}{\ensuremath{\mathrm{U}}}
\newcommand{\chara}{\ensuremath{\mathrm{char}}}
\newcommand{\ann}{\ensuremath{\mathsf{Ann}}}

\newcommand{\C}{\ensuremath{\mathcal{C}}}
\newcommand{\A}{\ensuremath{\mathcal{A}}}
\newcommand{\Ab}{\ensuremath{\mathsf{Ab}}}
\newcommand{\ab}{\ensuremath{\mathsf{ab}}}
\newcommand{\vect}{\ensuremath{\mathsf{vect}}}
\newcommand{\lie}{\ensuremath{(-)_\mathsf{Lie}}}
\newcommand{\XMod}{\ensuremath{\mathsf{XMod}}}
\newcommand{\Centr}{\ensuremath{\mathsf{Centr}}}
\newcommand{\xmod}{\ensuremath{(-)_{\mathsf{Peiff}}}}

\newcommand{\K}{\ensuremath{\mathbb{K}}}
\newcommand{\B}{\ensuremath{\mathcal{B}}}
\renewcommand{\b}{\ensuremath{\mathsf{b}}}
\newcommand{\cc}{\ensuremath{\mathsf{c}}}

\newcommand{\Z}{\ensuremath{\mathbb{Z}}}
\newcommand{\Gp}{\ensuremath{\mathsf{Gp}}}
\newcommand{\Mod}{\ensuremath{\mathsf{Mod}}}
\newcommand{\Lie}{\ensuremath{\mathsf{Lie}}}
\newcommand{\Lieg}{\ensuremath{\mathfrak{g}}}
\newcommand{\Liec}{\ensuremath{\mathfrak{c}}}
\newcommand{\Lieb}{\ensuremath{\mathfrak{b}}}
\newcommand{\Liea}{\ensuremath{\mathfrak{a}}}
\newcommand{\Leibniz}{\ensuremath{\mathsf{Leib}}}
\newcommand{\Vect}{\ensuremath{\mathsf{Vect}}}
\newcommand{\PXMod}{\ensuremath{\mathsf{PXMod}}}
\renewcommand{\hom}{\ensuremath{\mathrm{Hom}}}

\newcommand{\pr}{\ensuremath{\mathrm{pr}}}
\newcommand{\CExt}{\ensuremath{\mathsf{CExt}}}
\newcommand{\NAAlg}{\ensuremath{\mathsf{NAAlg}}}
\newcommand{\Ext}{\ensuremath{\mathsf{Ext}}}
\newcommand{\Arr}{\ensuremath{\mathsf{Arr}}}
\newcommand{\tensor}{\ensuremath{\otimes}}

\newdir{>>}{{}*!/3.5pt/:(1,-.2)@^{>}*!/3.5pt/:(1,+.2)@_{>}*!/7pt/:(1,-.2)@^{>}*!/7pt/:(1,+.2)@_{>}}
\newdir{ >>}{{}*!/8pt/@{|}*!/3.5pt/:(1,-.2)@^{>}*!/3.5pt/:(1,+.2)@_{>}}
\newdir{ |>}{{}*!/-3.5pt/@{|}*!/-8pt/:(1,-.2)@^{>}*!/-8pt/:(1,+.2)@_{>}}
\newdir{ >}{{}*!/-8pt/@{>}}
\newdir{>}{{}*:(1,-.2)@^{>}*:(1,+.2)@_{>}}
\newdir{<}{{}*:(1,+.2)@^{<}*:(1,-.2)@_{<}}

\newbox\pullbackbox
\setbox\pullbackbox=\hbox{\xy 0;<1mm,0mm>: \POS(4,0)\ar@{-} (0,0) \ar@{-} (4,4)
\endxy}
\def\pullback{\copy\pullbackbox}

\newbox\pushoutbox
\setbox\pushoutbox=\hbox{\xy 0;<1mm,0mm>: \POS(0,4)\ar@{-} (0,0) \ar@{-} (4,4)
\endxy}

\begin{document}

\title[A relative theory of universal central extensions]{A relative theory of\\ universal central extensions}

\author{Jos\'e Manuel Casas}
\email{jmcasas@uvigo.es}
\address{Dpto.\ de Matem\'atica Aplicada I, Universidad de Vigo, Escola Enxe\~nar\'ia Forestal, Campus Universitario A Xunquiera, E-36005 Pontevedra, Spain}
\thanks{The first author's research was supported by Ministerio de Ciencia e Innovaci\'on (grant number MTM2009-14464-C02-02, includes European FEDER support) and by Xunta de Galicia (grant number Incite09 207 215 PR)}

\author{Tim Van~der Linden}
\email{tim.vanderlinden@uclouvain.be}
\address{Institut de recherche en math\'ematique et physique, Universit\'e catholique de Louvain, chemin du cyclotron~2 bte~L7.01.02, B-1348 Louvain-la-Neuve, Belgium}
\address{CMUC, University of Coimbra, 3001--454 Coimbra, Portugal}
\thanks{The second author works as \emph{charg\'e de recherches} for Fonds de la Recherche Scientifique--FNRS. His research was supported by Centro de Matem\'atica da Universidade de Coimbra and by Funda\c c\~ao para a Ci\^encia e a Tecnologia (grant number SFRH/BPD/38797/2007).}

\begin{abstract}
Basing ourselves on Janelidze and Kelly's general notion of central extension, we study universal central extensions in the context of semi-abelian categories. Thus we unify classical, recent and new results in one conceptual framework. The theory we develop is relative to a chosen Birkhoff subcategory of the category considered: for instance, we consider groups vs.\ abelian groups, Lie algebras vs.\ vector spaces, precrossed modules vs.\ crossed modules and Leibniz algebras vs.\ Lie algebras. We consider a fundamental condition on composition of central extensions and give examples of categories which do, or do not, satisfy this condition.
\end{abstract}




\date{\today}

\maketitle

\section*{Introduction}
This article provides a unified framework for the study of universal central extensions, using techniques from categorical Galois theory~\cite{Janelidze:Pure} and, in particular, Janelidze and Kelly's relative theory of central extensions~\cite{Janelidze-Kelly}. Its aim is to make explicit the underlying unity of results in the literature (for groups, Leibniz algebras, precrossed modules, etc.\ \cite{ACL, Casas-Corral, Casas-Ladra, Cheng-Su, Gn, Gran-VdL}) and to unite them in a single, general setting. Thus a basic theory of universal central extensions is developed for all these special cases at once.

We work in a pointed Barr exact Goursat category $\A$ with a chosen Birkhoff subcategory $\B$ of $\A$; the universal central extensions of $\A$ are defined relative to the chosen $\B$. This is the minimal setting in which the theory of central extensions from~\cite{Janelidze-Kelly} can be used to obtain meaningful results on the relations between perfect objects and universal central extensions. (Indeed, we need $\A$ to be Barr exact and Goursat for the concepts of normal and central extension to coincide, and for split epimorphic central extensions to be trivial; and perfect objects can only be properly considered in a pointed context.)

The simultaneously \emph{categorical} and \emph{Galois theoretic} approach due to Janelidze and Kelly is based on, and generalises, the work of the Fr\"ohlich school~\cite{Froehlich, Furtado-Coelho, Lue} which focused on varieties of $\Omega$-groups. Recall~\cite{Higgins} that a \defn{variety of $\Omega$-groups} is a variety of universal algebras which has amongst its operations and identities those of the variety of groups but has just one constant; furthermore, a Birkhoff subcategory of a variety is the same thing as a subvariety.

In order to construct universal central extensions, we further narrow the context to that of semi-abelian categories with enough projectives~\cite{Borceux-Bourn, Janelidze-Marki-Tholen}, which still includes all varieties of $\Omega$-groups. We need a good notion of short exact sequence to construct the centralisation of an extension, and the existence of projective objects gives us weakly universal central extensions. The switch to semi-abelian categories also allows us to make the connection with existing homology theories~\cite{EverHopf, EverVdL1, EverVdL2} and to prove some classical recognition results for universal central extensions.

Although some examples (for instance groups vs.\ abelian groups and Lie algebras vs.\ vector spaces) are \emph{absolute}, meaning that they fit into the theory relative to the subcategory $\Ab(\A)$ of all abelian objects, others are not: precrossed modules vs.\ crossed modules, and Leibniz algebras vs.\ Lie algebras, for instance. In this absolute case, some results were already investigated in~\cite{Gran-VdL}; they appear as special cases of our general theory. 

What surprised us most of all while developing this relative theory of universal central extensions is that the semi-abelian context turns out to be too weak for some of the most basic results, valid for classical examples such as groups and Lie algebras, to hold. We have to impose an additional condition which we chose to call \emph{\UCEspecial}. It asks that if $B$ is an $\ab$-perfect object in the category $\A$ and $f\colon{B\to A}$ and $g\colon{C\to B}$ are $\ab$-central extensions, then also the composite extension $f\comp g$ is $\ab$-central. Under \UCE, standard recognition results such as Theorem~\ref{Theorem-Homology} hold true. But \UCE\ is also necessary. This immediately gives rise to the following question, which is not yet fully answered in this paper: \emph{When does \UCE\ hold?} We can give examples and counterexamples, but thus far there is no elementary characterisation. Thus we propose the following problem: \emph{To find good minimal hypotheses for \UCE}.

The text is structured as follows. In the first section we develop that part of the theory which does not depend on the existence of either projective objects or short exact sequences. Here we work in pointed Barr exact Goursat categories. We sketch the context and recall the basic definitions of \emph{perfect object} and \emph{(universal) central extension}. Some of the simpler correspondences between them are developed, as for instance Proposition~\ref{Proposition-Universal-vs-Weak-Universal} on the universality of a central extension vs.\ perfectness of its domain. Further results are obtained in the setting of semi-abelian categories with enough projectives. In Section~\ref{Section-UCE-Construction} we prove that any perfect object admits a universal central extension (Theorem~\ref{Theorem-Universal-Central-Extension}). In Section~\ref{Section-Nested-Birkhoff-Subcategories} we consider the case of nested Birkhoff subcategories $\C\subset\B\subset\A$. Given a perfect object $B$ of~$\B$ we obtain a short exact sequence comparing the second homology of~$B$, viewed as an object of~$\B$, with the second homology of $B$, viewed as an object of $\A$ (Proposition~\ref{Proposition-Comparison}). In Section~\ref{Section-Recognition} we show that, when $\Ab(\A)\subset \B$ and \UCE\ holds, a central extension is universal exactly when its domain is perfect and projective with respect to all central extensions (Proposition~\ref{Proposition-Universal-vs-Weak-Universal'}), and we also refine the connections with semi-abelian homology obtained in Section~\ref{Section-UCE-Construction} (Theorem~\ref{Theorem-Homology}). In Section~\ref{Section-Examples} we show how the theory unifies existing with new results by explaining the examples of groups vs.\ abelian groups, Leibniz algebras vs.\ Lie algebras vs.\ vector spaces, and precrossed modules vs.\ crossed modules vs.\ abelian crossed modules. We also give counterexamples related to \UCE. Finally, in Section~\ref{Section-Final-Remarks}, we end with some additional remarks: on possible characterisations of \UCE, on stem extensions, and on universal higher central extensions.

\section{Basic definitions and results}\label{Section-Basic-Definitions}
In their article~\cite{Janelidze-Kelly}, Janelidze and Kelly introduced a general theory of relative central extensions in the context of exact Goursat categories. This is the theory we shall be considering here, focusing on the induced relative notion of universal central extension. We give an overview of the needed definitions and prove some preliminary results on the relation between universal central extensions and perfect objects. In the following section we shall narrow the context to semi-abelian categories with enough projectives in order to prove the existence of universal central extensions.

\subsection{Barr exact Goursat categories}
Recall that a \defn{regular epimorphism} is a coequaliser of some pair of arrows. A category is called \defn{regular} when it is finitely complete with coequalisers of kernel pairs and with pullback-stable regular epimorphisms. In a regular category, any morphism may be factored as a regular epimorphism followed by a monomorphism, and this \defn{image factorisation} is unique up to isomorphism. A category is said to be \defn{Barr exact} if and only if it is regular and such that any internal equivalence relation is a kernel pair. 

Next to Barr exactness, the theory of central extensions considered in~\cite{Janelidze-Kelly} needs the surrounding category to satisfy the \emph{Goursat property}. (Then, for instance, the concepts of normal extension and central extension coincide, and every split epimorphic central extension is trivial. Both of these facts are crucial in what follows.) A Barr exact category is called \defn{Goursat} when for every pair of equivalence relations $R$, $S$ on an object $X$ the condition $SRS=RSR$ holds. For most examples the slightly less general and better known context of exact \defn{Mal'tsev} categories suffices: here any internal reflexive relation is an equivalence relation or, equivalently, the condition $SR=RS$ holds for all equivalence relations $R$, $S$ on an object $X$. A variety is Mal'tsev in this sense if and only if it is a \emph{Mal'tsev variety}. Moreover, a Barr exact category is Mal'tsev if and only if the pushout of a regular epimorphism along a regular epimorphism always exists, and the comparison morphism to the induced pullback is also a regular epimorphism~\cite{Carboni-Kelly-Pedicchio}. See~\cite{Janelidze-Kelly} for further details.

\subsection{Birkhoff subcategories}
The notion of central extension introduced in~\cite{Janelidze-Kelly} is \emph{relative}, being defined with respect to a chosen subcategory $\B$ of the category $\A$ considered.

Let $\A$ be a Barr exact Goursat category. A \defn{Birkhoff subcategory} $\B$ of~$\A$ is a full and reflective subcategory which is closed under subobjects and regular quotients. We write the induced adjunction as
\begin{equation}\label{Birkhoff-Adjunction}
\xymatrix@1{{\A} \ar@<1ex>[r]^-{\b} \ar@{}[r]|-{\perp} & {\B} \ar@<1ex>[l]^-{\supset}}
\end{equation}
and denote its unit by $\eta\colon{1_{\A}\To \b}$. A Birkhoff subcategory of a variety of universal algebras is the same thing as a subvariety. If $\A$ is finitely complete Barr exact Goursat then so is any Birkhoff subcategory $\B$ of $\A$.

For a given full, replete and reflective subcategory $\B$, being closed under subobjects is equivalent to the components $\eta_{A}$ of the unit of the adjunction being regular epimorphisms. If now $\B$ is full, reflective and closed under subobjects then the Birkhoff property of $\B$ (that is, closure under quotients) is equivalent to the following condition: given any regular epimorphism $f\colon{B\to A}$ in $\A$,
the induced square of regular epimorphisms
\begin{equation}\label{Birkhoff-Square}
\vcenter{\xymatrix{B \ar[r]^-{f} \ar[d]_-{\eta_{B}} & A \ar[d]^-{\eta_{A}}\\
\b(B) \ar[r]_-{\b(f)} & \b(A)}}
\end{equation}
is a pushout.

 From now on, $\B$ will be a fixed Birkhoff subcategory of a chosen Barr exact Goursat category $\A$.

\subsection{Extensions and central extensions}\label{Subsection-Central-Extensions}

An \defn{extension} in $\A$ is a regular epimorphism. A morphism of extensions is a commutative square between them, and thus we obtain the category $\Ext(\A)$ of extensions in $\A$. 

With respect to the Birkhoff subcategory $\B$, there are notions of \emph{trivial}, \emph{normal} and \emph{central} extension. An extension $f\colon{B\to A}$ in $\A$ is said to be \defn{trivial (with respect to~$\B$)} or \defn{$\b$-trivial} if and only if the induced square~\eqref{Birkhoff-Square} is a pullback. The extension $f$ is \defn{normal (with respect to $\B$)} or \defn{$\b$-normal} if and only if one of the projections $f_{0}$, $f_{1}$ in the kernel pair $(B\times_A B,f_{0},f_{1})$ of $f$ is $\b$-trivial. That is to say, $f$ is normal with respect to $\B$ precisely when in the diagram
\[
\xymatrix{B\times_A B \ar@<.5ex>[r]^-{f_{0}} \ar@<-.5ex>[r]_-{f_{1}} \ar[d]_-{\eta_{B\times_A B}} & B \ar[r]^-{f} \ar[d]^-{\eta_{B}} & A\\
\b(B\times_A B) \ar@<.5ex>[r]^-{\b(f_{0})} \ar@<-.5ex>[r]_-{\b(f_{1})} & \b(B)}
\]
both commutative squares are pullbacks, since one of them being a pullback implies that so is the other. Finally, $f$ is said to be \defn{central (with respect to~$\B$)} or \defn{$\b$-central} if and only if there exists an extension $g\colon{C\to A}$ such that the pullback~$g^{*}(f)$ of $f$ along $g$ is~$\b$-trivial.

\begin{remark}\label{Remark-Central-and-Trivial}
Clearly, every normal extension is central; in the present context, the converse also holds, and thus the concepts of normal and central extension coincide. Furthermore, a split epimorphism is a trivial extension if and only if it is a central extension~\cite[Theorem~4.8]{Janelidze-Kelly}. Finally, central extensions are pullback-stable~\cite[Proposition~4.3]{Janelidze-Kelly}.
\end{remark}

\begin{remark}\label{Remark-Galois-Structure}
Together with the classes $|\Ext(\A)|$ and $|\Ext(\B)|$ of extensions in $\A$ and in~$\B$, the adjunction~\eqref{Birkhoff-Adjunction} forms a \emph{Galois structure}~\cite{Janelidze:Pure}
\[
\Gamma=\bigl(\xymatrix@1{{\A} \ar@<1ex>[r]^-{\b} \ar@{}[r]|-{\perp} & {\B} \ar@<1ex>[l]^-{\supset}},|\Ext(\A)|,|\Ext(\B)|\bigr).
\]
\end{remark}

\subsection{Pointed categories}

In what follows it will be crucial that the terminal object~$1$ of $\A$ is also initial, so that the category~$\A$ is \defn{pointed}. When this happens, the object $1=0$ is called the \defn{zero object} of $\A$. A morphism $f$ is \defn{zero} when it factors over the zero object.

Since the reflector $\b$ always preserves pullbacks of split epimorphisms along split epimorphisms, in the pointed case it also preserves products.

 From now on, $\A$ will be a fixed pointed exact Goursat category; any Birkhoff subcategory $\B$ of $\A$ is also pointed exact Goursat.

\subsection{Perfect objects}

An object $P$ of $\A$ is called \defn{perfect (with respect to $\B$)} or \defn{$\b$-perfect} when $\b(P)$ is the zero object $0$ of $\B$. If $f\colon{B\to A}$ is an extension and~$B$ is $\b$-perfect then so is $A$, because the reflector $\b$ preserves regular epimorphisms, and a regular quotient of zero is zero.

\begin{lemma}\label{Lemma-Perfect-and-Central}
Let $P$ be a $\b$-perfect object and $f\colon{B\to A}$ an extension.
\begin{enumerate}
\item If $f$ is $\b$-trivial then the map
\[
\hom(P,f)=f\comp(-)\colon{\hom(P,B)\to \hom(P,A)}
\]
is a bijection;
\item if $f$ is $\b$-central then $\hom(P,f)$ is an injection.
\end{enumerate}
If $\hom(P,f)$ is an injection for every $\b$-trivial extension $f$ then $P$ is $\b$-perfect.
\end{lemma}
\begin{proof}
The extension $f$ being $\b$-trivial means that the square~\eqref{Birkhoff-Square} is a pullback. If $b_{0}$, $b_{1}\colon{P\to B}$ are morphisms such that $f\comp b_{0}=f\comp b_{1}$ then $b_{0}$ is equal to $b_{1}$ by the uniqueness in the universal property of this pullback: indeed also $\eta_{B}\comp b_{0}=\b(b_{0})\comp\eta_{P}=0=\b(b_{1})\comp\eta_{P}=\eta_{B}\comp b_{1}$. Thus we see that $\hom(P,f)$ is injective.
This map is also surjective, since any morphism $a\colon{P\to A}$ is such that $\eta_{A}\comp a=\b(a)\comp \eta_{P}=0=\b(f)\comp 0$ and thus induces a morphism $b\colon{P\to B}$ for which $f\comp b=a$.

Statement 2 follows from 1 because the functor $\hom(P,-)$ preserves kernel pairs, and a map is an injection if and only if its kernel pair projections are bijections.

As to the converse: the morphism $!_{\b(P)}\colon{\b(P)\to 0}$ is a $\b$-trivial extension; since $!_{\b(P)}\comp \eta_{P}=0=!_{\b(P)}\comp 0$, the assumption implies that $\eta_{P}$ is zero, which means that $P$ is $\b$-perfect.
\end{proof}

\subsection{Universal central extensions}
For any object $A$ of $\A$, let $\Centr_{\b}(A)$ denote the category of all $\b$-central extensions of $A$: the full subcategory of the slice category $(\A\downarrow A)$ determined by the central extensions. A (weakly) initial object of this category $\Centr_{\b}(A)$ is called a \defn{(weakly) universal $\b$-central extension} of~$A$. A $\b$-central extension $u\colon{U\to A}$ is weakly universal when for every $\b$-central extension $f\colon{B\to A}$ there exists a morphism $\overline{f}$ from $u$ to $f$, that is, such that $f\comp \overline{f}=u$. Furthermore, $u$ is universal when this induced morphism $\overline{f}$ is unique. Note also that, up to isomorphism, an object admits at most one universal $\b$-central extension.

\begin{lemma}\label{UCE-then-Perfect}
If $u\colon{U\to A}$ is a universal $\b$-central extension then the objects $U$ and $A$ are $\b$-perfect.
\end{lemma}
\begin{proof}
Since the first projection $\pr_{A}\colon {A\times \b(U)\to A}$ is a trivial extension, by Remark~\ref{Remark-Central-and-Trivial} it is central. By the hypothesis that $u$ is universal, there exists just one morphism $\langle u,v\rangle\colon {U\to A\times \b(U)}$ such that $\pr_{A}\comp\langle u,v\rangle=u$. But then $0\colon {U\to \b(U)}$ is equal to $\eta_{U}\colon {U\to \b(U)}$, and $\b(U)=0$. Since a regular quotient of a perfect object is perfect, this implies that both $U$ and $A$ are $\b$-perfect.
\end{proof}

\begin{proposition}\label{Proposition-Universal-vs-Weak-Universal}
Let $\A$ be a pointed Barr exact Goursat category and~$\B$ a Birkhoff subcategory of $\A$. Let $u\colon{U\to A}$ be a $\b$-central extension. Between the following conditions, the implications $1\Leftrightarrow 2 \Leftrightarrow 3\Rightarrow 4\Leftrightarrow 5$ hold:
\begin{enumerate}
\item $U$ is $\b$-perfect and every $\b$-central extension of $U$ splits;
\item $U$ is $\b$-perfect and projective with respect to all $\b$-central extensions;
\item for every $\b$-central extension $f\colon{B\to A}$, the map
\[
\hom(U,f)\colon{\hom(U,B)\to \hom(U,A)}
\]
is a bijection; 
\item $U$ is $\b$-perfect and $u$ is a weakly universal $\b$-central extension;
\item $u$ is a universal $\b$-central extension.
\end{enumerate}
\end{proposition}
\begin{proof}
Suppose that 1 holds. To prove 2, let $f\colon{B\to A}$ be a $\b$-central extension and $g\colon{U\to A}$ a morphism. Then the pullback $g^{*}(f)\colon{\overline{B}\to U}$ of $f$ along $g$ is still $\b$-central; hence $g^{*}(f)$ admits a splitting $s\colon{U\to \overline{B}}$, and $(f^{*}(g))\comp s$ is the required morphism ${g\to f}$. Conversely, given a $\b$-central extension $f\colon{B\to U}$, the projectivity of~$U$ yields a morphism $s\colon{U\to B}$ such that $f\comp s=1_{U}$. 

Conditions 2 and 3 are equivalent by Lemma~\ref{Lemma-Perfect-and-Central}.

Condition 3 implies condition 5: given a $\b$-central extension $f\colon{B\to A}$ of $A$, there exists a unique morphism $\overline{f}\colon{U\to B}$ that satisfies $f\comp \overline{f}=u$.

Finally, 4 and 5 are equivalent by Lemma~\ref{Lemma-Perfect-and-Central} and Lemma~\ref{UCE-then-Perfect}.
\end{proof}

\begin{remark}
To prove that condition 4 implies 3 we would require $U$ itself to admit a universal $\b$-central extension, which need not be the case in the present context. In fact, even if such a universal $\b$-central extension of $U$ does exist then the above five conditions may or may not be equivalent: see Section~\ref{Section-Recognition}.
\end{remark}

\section{The universal central extension construction}\label{Section-UCE-Construction}
Our aim is now to prove that every perfect object admits a universal central extension. To do so, a richer categorical context is needed; for instance, a good notion of short exact sequence will be crucial in the construction of the centralisation of an extension and in the passage to a perfect subobject of an object. The existence of projective objects will also become important now: they will give us weakly universal central extensions. We switch to the framework of semi-abelian categories with enough projectives.

\subsection{Semi-abelian categories}
A pointed and regular category is called \defn{Bourn protomodular} when the \defn{(Regular) Short Five Lemma} holds: this means that for any commutative diagram
\begin{equation}\label{Short-Five-Lemma}
\vcenter{\xymatrix{\Ker(f') \ar[r]^-{\ker f'} \ar[d]_-k & B' \ar[r]^-{f'} \ar[d]_-b & A' \ar[d]^-a \\ \Ker(f) \ar[r]_-{\ker f} & B \ar[r]_-{f} & A}}
\end{equation}
such that $f$ and $f'$ are regular epimorphisms, $k$ and $a$ being isomorphisms implies that $b$ is an isomorphism. A category is \defn{semi-abelian} when it is pointed, Barr exact and Bourn protomodular with binary coproducts~\cite{Janelidze-Marki-Tholen}. A variety of $\Omega$-groups is always a semi-abelian category. A semi-abelian category is always a Mal'tsev category (hence it is also Goursat)~\cite{Borceux-Bourn}.

Since a regular epimorphism is always the cokernel of its kernel in a semi-abelian category, an appropriate notion of short exact sequence exists. A \defn{short exact sequence} is any sequence
\[
\xymatrix{K \ar[r]^-{k} & B \ar[r]^-{f} & A}
\]
that satisfies $k=\ker f$ and $f=\coker k$. We denote this situation by
\begin{equation}\label{Short-Exact-Sequence}
\xymatrix{0 \ar[r] & K \ar@{{ |>}->}[r]^-{k} & B \ar@{ >>}[r]^-{f} & A \ar[r] & 0.}
\end{equation}

\begin{lemma}\cite{Bourn1991, Bourn2001}\label{Lemma-Pullback}
Consider a morphism of short exact sequences~\eqref{Short-Five-Lemma}.
\begin{enumerate}
\item The right hand side square $f\comp b=a\comp f'$ is a pullback iff $k$ is an isomorphism.
\item The left hand side square $\ker f\comp k=b\comp \ker f'$ is a pullback iff $a$ is mono.\hfill\qed
\end{enumerate}
\end{lemma}

The first statement implies that any pullback square between regular epimorphisms (that is, any square $f\comp b=a\comp f'$ as in~\eqref{Short-Five-Lemma}) is a pushout. It is also well known that the regular image of a kernel is a kernel~\cite{Janelidze-Marki-Tholen}. In any semi-abelian category, classical homological lemma's such as the Snake Lemma and the $3\times3$ Lemma are valid; for further details and many other results we refer the reader to~\cite{Janelidze-Marki-Tholen, Borceux-Bourn}.

As of now $\A$ will be a semi-abelian category and $\B$ a Birkhoff subcategory of $\A$.

\subsection{Commutators and centralisation}\label{Subsection-Commutators-and-Centralisation}
The kernel $\mu$ of the unit $\eta$ of the adjunction~\eqref{Birkhoff-Adjunction} gives rise to a ``zero-dimensional'' commutator: for any object~$A$ of~$\A$,
\[
\xymatrix{0 \ar[r] & [A,A]_{\b} \ar@{{ |>}->}[r]^-{\mu_{A}} & A \ar@{ >>}[r]^-{\eta_{A}} & \b(A) \ar[r] & 0}
\]
is a short exact sequence in $\A$; hence $A$ is an object of $\B$ if and only if $[A,A]_{\b}=0$. On the other hand, an object~$A$ of~$\A$ is $\b$-perfect precisely when $[A,A]_{\b}=A$. This construction defines a functor $[-,-]_{\b}\colon{\A\to \A}$ and a natural transformation~$\mu\colon{[-,-]_{\b}\To 1_{\A}}$. The functor $[-,-]_{\b}$ preserves regular epimorphisms; we recall the argument. Given a regular epimorphism $f\colon{B\to A}$, by the Birkhoff property, the induced square of regular epimorphisms on the right
\[
\xymatrix{0 \ar[r] & [B,B]_{\b} \ar@{{ |>}->}[r]^-{\mu_{B}} \ar@{.>}[d]_-{[f,f]_{\b}} & B \ar@{ >>}[d]_-{f} \ar@{ >>}[r]^-{\eta_{B}} & \b(B) \ar@{ >>}[d]^-{\b(f)} \ar[r] & 0\\
0 \ar[r] & [A,A]_{\b} \ar@{{ |>}->}[r]_-{\mu_{A}} & A \ar@{ >>}[r]_-{\eta_{A}} & \b(A) \ar[r] & 0}
\]
is a pushout---but this is equivalent to $[f,f]_{\b}$ being a regular epimorphism~\cite[Corollary 5.7]{EverVdL1}.

Lemma~\ref{Lemma-Pullback} implies that an extension $f$ as in~\eqref{Short-Exact-Sequence} is $\b$-central if and only if either one of the morphisms $[f_{0},f_{0}]_{\b}$, $[f_{1},f_{1}]_{\b}$ is an isomorphism, which happens exactly when they coincide, $[f_{0},f_{0}]_{\b}=[f_{1},f_{1}]_{\b}$.
\[
\xymatrix{& [K,B]_{\b} \ar@{{ |>}->}[d]_-{\ker [f_{0},f_{0}]_{\b}} \ar@{{ |>}.>}[ddr] \\
0 \ar[r] & [B\times_A B,B\times_A B]_{\b} \ar@{{ |>}->}[r]^-{\mu_{B\times_A B}} \ar@<-.5ex>@{ >>}[d]_-{[f_{0},f_{0}]_{\b}} \ar@<.5ex>@{ >>}[d]^-{[f_{1},f_{1}]_{\b}} & B\times_A B \ar@<-.5ex>@{ >>}[d]_-{f_{0}} \ar@<.5ex>@{ >>}[d]^-{f_{1}} \ar@{ >>}[r]^-{\eta_{B\times_A B}} & \b(B\times_A B) \ar@<-.5ex>@{ >>}[d]_-{\b(f)_{0}} \ar@<.5ex>@{ >>}[d]^-{\b(f)_{1}} \ar[r] & 0\\
0 \ar[r] & [B,B]_{\b} \ar@{{ |>}->}[r]_-{\mu_{B}} & B \ar@{ >>}[r]_-{\eta_{B}} & \b(B) \ar[r] & 0}
\]
Hence the kernel $[K,B]_{\b}$ of $[f_{0},f_{0}]_{\b}$ measures how far $f$ is from being central: indeed, $f$ is $\b$-central if and only if $[K,B]_{\b}$ is zero. 

\begin{remark}\label{Remark-Subobject}
This explains, for instance, why a sub-extension of a central extension is central. It is worth remarking here that a morphism of extensions $(b_{0},a_{0})$ as in~\eqref{Morphisms} below is a monomorphism if and only if such is $b_{0}$.
\end{remark}

The ``one-dimensional'' commutator $[K,B]_{\b}$ may be considered as a normal subobject of $B$ via the composite $\mu_{B}\comp [f_{1},f_{1}]_{\b}\comp \ker [f_{0},f_{0}]_{\b}\colon{[K,B]_{\b}\to B}$. Thus the Galois structure $\Gamma$ from Remark~\ref{Remark-Galois-Structure} induces a new adjunction
\[
\xymatrix@1{{\Ext(\A)} \ar@<1ex>[r]^-{\b_{1}} \ar@{}[r]|-{\perp} & {\CExt_{\b}(\A)}, \ar@<1ex>[l]^-{\supset}}
\]
where $\CExt_{\b}(\A)$ is the full reflective subcategory of $\Ext(\A)$ determined by the $\b$-central extensions. Given an extension $f\colon{B\to A}$ with kernel $K$, its \defn{centralisation} $\b_{1}(f)\colon{B/[K,B]_{\b}\to A}$ is obtained through the diagram with exact rows
\[
\xymatrix{0 \ar[r] & [K,B]_{\b} \ar@{{ |>}->}[r] \ar@{ >>}[d] & B \ar@{ >>}[r] \ar@{ >>}[d]_-{f} & \tfrac{B}{[K,B]_{\b}} \ar[r] \ar@{ >>}[d]^-{\b_{1}(f)} & 0\\
& 0 \ar@{{ |>}->}[r] & A \ar@{=}[r] & A \ar[r] & 0.}
\]
Considering this diagram as a short exact sequence
\[
\xymatrix{0 \ar[r] & \Ker(\eta^{1}_{f}) \ar@{{ |>}->}[r]^-{\mu^{1}_{f}} & f \ar@{ >>}[r]^-{\eta^{1}_{f}} & \b_{1}(f) \ar[r] & 0}
\]
in the semi-abelian category of arrows $\Arr(\A)$---morphisms here are commutative squares---we obtain a description of the unit $\eta^{1}$ of the adjunction and its kernel~$\mu^{1}$.

\subsection{Baer invariants}\label{Subsection-Baer-Invariants}
We recall the basic definitions of the theory of Baer invariants~\cite[Definition~3.1 and 3.3]{EverVdL1}. Two given morphisms of extensions 
\begin{equation}\label{Morphisms}
\vcenter{\xymatrix{B' \ar@{-{ >>}}[d]_-{f'} \ar@<-0.5 ex>[r]_-{b_{1}}\ar@<0.5 ex>[r]^-{b_{0}} & B \ar@{-{ >>}}[d]^-{f}\\
A' \ar@<-0.5 ex>[r]_-{a_{1}}\ar@<0.5 ex>[r]^-{a_{0}} & A}}
\end{equation}
$(b_{0},a_{0})$ and $(b_{1},a_{1})\colon {f'\to f}$ are said to be \defn{homotopic} when $a_{0}=a_{1}$. A \defn{Baer invariant} is a functor $F\colon {\Ext(\A) \to \A} $ which makes homotopic morphisms of extensions equal: if $(b_{0},a_{0})\simeq (b_{1},a_{1})$ then $F(b_{0},a_{0})=F(b_{1},a_{1})$. Such a functor sends homotopically equivalent extensions to isomorphic objects. 

For instance, the functor ${\Ext(\A) \to \A}$ that maps an extension
\begin{equation}\label{Extension}
\xymatrix{0 \ar[r] & K \ar@{{ |>}->}[r]^{k} & B \ar@{-{ >>}}[r]^{f} & A \ar[r] & 0}
\end{equation}
to the quotient $[B,B]_{\b}/[K,B]_{\b}$ is an example of a Baer invariant, as is the functor which maps this extension to the quotient ${(K\meet [B,B]_{\b})}/{[K,B]_{\b}}$. See~\cite{EverVdL1}, and in particular its Proposition~4.6, for further details.

\subsection{Existence of a weakly universal central extension}\label{Subsection-Enough-Projectives}
We say that $\A$ \defn{has weakly universal central extensions} (for some Birkhoff subcategory $\B$ of $\A$) when every object of $\A$ admits a weakly universal $\b$-central extension. This happens, for instance, when $\A$ has enough (regular) projectives, so that for any object~$A$ of~$\A$, there exists a regular epimorphism $f\colon{B\to A}$ with $B$ projective, a \defn{(projective) presentation} of $A$.

\begin{lemma}\label{Lemma-Weakly-UCE}
If the category $\A$ is semi-abelian with enough projectives then it has weakly universal central extensions for any Birkhoff subcategory $\B$.
\end{lemma}
\begin{proof}
Given an object $A$ of $\A$, the category $\Centr_{\b}(A)$ has a weakly initial object: given a projective presentation $f\colon{B\to A}$ with kernel $K$, its centralisation $\b_{1}(f)\colon{B/[K,B]_{\b}\to A}$ is weakly initial. Indeed, any $\b$-central extension $g\colon {C\to A}$ induces a morphism ${\b_{1}(f)\to g}$ in $\Centr_{\b}(A)$, as the object $B$ is projective.
\end{proof}

\subsection{The Schur multiplier}\label{Subsection-Schur-Multiplier}
Let $A$ be an object of $\A$ and $f\colon{B\to A}$ a projective presentation with kernel~$K$. The induced objects
\[
\frac{[B,B]_{\b}}{[K,B]_{\b}}\qquad\qquad \text{and} \qquad\qquad \frac{K\meet [B,B]_{\b}}{[K,B]_{\b}}
\]
are independent of the chosen projective presentation of $A$ as explained above. Hence the following makes sense:

\begin{definition}\label{Definition-Homology}
By analogy with classical homology theories, the latter object ${(K\meet [B,B]_{\b})}/{[K,B]_{\b}}$ is called the \defn{second homology object} or the \defn{Schur multiplier} of $A$ \defn{(relative to $\B$)} and is written $\H_{2}(A,\b)$. We write $\U(A,\b)$ for the object ${[B,B]_{\b}}/{[K,B]_{\b}}$, and $\H_{1}(A,\b)$ will denote the reflection $\b(A)$ of $A$ into~$\B$.
\end{definition}

\begin{remark}\label{Remark-Homology}
The objects $\H_{2}(A,\b)$ and $\H_{1}(A,\b)$ are genuine homology objects: if~$\A$ is a semi-abelian monadic category then they may be computed using comonadic homology as in~\cite{EverVdL2}---note that the monadicity here implies existence of enough projectives. In any case, they fit into the homology theory worked out in~\cite{EverHopf}. By Theorem~5.9 in~\cite{EverVdL1}, any short exact sequence~\eqref{Short-Exact-Sequence} induces a five-term exact sequence
\[
\xymatrix{\H_{2}(B,\b) \ar[r] & \H_{2}(A,\b) \ar[r] & \frac{K}{[K,B]_{\b}} \ar[r] & \H_{1}(B,\b) \ar@{-{ >>}}[r] & \H_{1}(A,\b) \ar[r] & 0,}
\]
the tail of a long exact homology sequence~\cite{EverHopf}. This is a relative generalisation of the Stallings--Stammbach sequence for groups, which is recovered when $\b=\ab\colon{\A=\Gp\to \B=\Ab}$ is the abelianisation functor---in fact, it is a categorical version of the analogous results which may be found in~\cite{Froehlich, Furtado-Coelho, Lue}.
\end{remark}

\subsection{Existence of a universal central extension}\label{UCE}
The Baer invariants from~\ref{Subsection-Baer-Invariants} may now be considered for all weakly universal $\b$-central extensions of an object $A$: indeed, any two such extensions of $A$ are always homotopically equivalent. Since for any weakly universal $\b$-central extension~\eqref{Extension} the commutator $[K,B]_{\b}$ is zero, the objects
\[
[B,B]_{\b}\qquad\qquad \text{and} \qquad\qquad K\meet [B,B]_{\b}
\]
are independent of the chosen weakly universal central extension of $A$. (Here, as in~\cite{Janelidze:Hopf}, the Hopf formula becomes $\H_{2}(A,\b)=K\meet [B,B]_{\b}$. Also note that $\U(A,\b)=[B,B]_{\b}$.)

We are now ready to prove that, if $A$ is $\b$-perfect, then a universal $\b$-central extension of $A$ does exist. This is a relative version of Proposition~4.1 in~\cite{Gran-VdL}.

\subsection{The perfect subobject}\label{Perfect-Subobject}
When there are weakly universal central extensions, any central extension of a perfect object contains a subobject with a perfect domain. We prove this in two steps: first for weakly universal central extensions, then in general. This implies that any perfect object admits a universal central extension when weakly universal central extensions exist.

\begin{lemma}\label{Lemma-Perfect-Subobject}
Suppose $\A$ is a semi-abelian category with a Birkhoff subcategory~$\B$.
Then any weakly universal $\b$-central extension of a $\b$-perfect object contains a subobject with a $\b$-perfect domain.
\end{lemma}
\begin{proof}
Let~\eqref{Extension} be a weakly universal $\b$-central extension of a $\b$-perfect object $A$. Since $\mu_{A}$ is an isomorphism and $[f,f]_{\b}$ is a regular epimorphism, the morphism $f\comp\mu_{B}=\mu_{A}\comp [f,f]_{\b}$ in the induced diagram with exact rows
\[
\xymatrix{
0 \ar[r] & K\meet [B,B]_{\b} \ar@{}[rd]|<<<{\pullback} \ar@{{ |>}->}[r] \ar@{{ |>}->}[d] & [B,B]_{\b} \ar@{{ |>}->}[d]^-{\mu_{B}} \ar@{-{ >>}}[r]^-{f\circ \mu_{B}} & A \ar[r] \ar@{=}[d] & 0\\
0\ar[r] & {K} \ar@{{ |>}->}[r] & {B} \ar@{ >>}[r]_-{f} & A \ar[r] & 0}
\]
is also a regular epimorphism. The extension $f\comp \mu_{B}$ is $\b$-central as a subobject of the $\b$-central extension $f$; its weak universality is clear. By Proposition~\ref{Proposition-Universal-vs-Weak-Universal}, the object $[B,B]_{\b}$ is $\b$-perfect, because the extensions $f\comp \mu_{B}$ and $f$ are homotopically equivalent, so that $[B,B]_{\b}\cong[[B,B]_{\b},[B,B]_{\b}]_{\b}$.
\end{proof}

\begin{lemma}\label{Lemma-Commutator-Perfect}
Let $\A$ be a semi-abelian category with weakly universal central extensions for a Birkhoff subcategory $\B$ of $\A$. If $f\colon {B\to A}$ is a $\b$-central extension of a $\b$-perfect object~$A$, then $[B,B]_{\b}$ is also $\b$-perfect.
\end{lemma}
\begin{proof}
The object $B$ admits a weakly universal $\b$-central extension $v\colon {V\to B}$; then the centralisation $w\colon{W\to A}$ of the resulting composite $f\comp v$ is a weakly universal $\b$-central extension. Indeed, given any $\b$-central extension $g\colon {C\to A}$, there is a factorisation $\overline{f^{*}(g)}\colon{V\to B\times_C A}$ of $v$ through the pullback $f^{*}(g)\colon{B\times_C A\to B}$ of $g$ along~$f$, and then the composite $(g^{*}(f))\comp (\overline{f^{*}(g)})\colon {V\to C}$ yields the needed morphism ${w\to g}$ by the universal property of the centralisation functor.

The comparison ${W\to B}$ is a regular epimorphism, hence so is the induced morphism ${[W,W]_{\b}\to [B,B]_{\b}}$; but a regular quotient of a perfect object is perfect.
\end{proof}

\begin{theorem}\label{Theorem-Universal-Central-Extension}
Let $\A$ be a semi-abelian category with enough projectives and $\B$ a Birkhoff subcategory of $\A$. An object $A$ of $\A$ is $\b$-perfect if and only if it admits a universal $\b$-central extension. Moreover, this universal $\b$-central extension may be chosen in such a way that it occurs in a short exact sequence
\[
\xymatrix{0 \ar[r] & \H_{2}(A,\b) \ar@{{ |>}->}[r] & \U(A,\b) \ar@{ >>}[r]^-{u^{\b}_{A}} & A \ar[r] & 0.}
\]
\end{theorem}
\begin{proof}
If an object admits a universal $\b$-central extension then it is $\b$-perfect by Lemma~\ref{UCE-then-Perfect}. Conversely, let~\eqref{Extension} be a weakly universal central extension of a $\b$-perfect object $A$ (Lemma~\ref{Lemma-Weakly-UCE}). Then by Lemma~\ref{Lemma-Perfect-Subobject} it admits a (weakly universal central) subobject with a $\b$-perfect domain. By Proposition~\ref{Proposition-Universal-vs-Weak-Universal}, this subobject is also universal. The shape of the short exact sequence follows from~\ref{UCE}.
\end{proof}

A variation on this proof allows us to show that the uniquely induced comparison morphism from a universal central extension to a central extension with a perfect domain is always a regular epimorphism.

\begin{proposition}\label{Proposition-Quotient-of-UCE}
Let $\A$ be a semi-abelian category with enough projectives and $\B$ a Birkhoff subcategory of $\A$. If $f\colon {B\to A}$ is a $\b$-central extension with a $\b$-perfect domain~$B$, then $f$ is a quotient of a universal $\b$-central extension.
\end{proposition}
\begin{proof}
The construction in the proof of Theorem~\ref{Theorem-Universal-Central-Extension} may be adapted to the given extension $f$ in such a way that the resulting morphism ${u\to f}$ is a regular epimorphism. We take a projective presentation $p\colon{P\to B}$ and use the composite $f\comp p\colon {P\to A}$ as a projective presentation of $A$. After centralisation we obtain a weakly universal $\b$-central extension $v\colon {V\to A}$ as in Lemma~\ref{Lemma-Weakly-UCE} and a regular epic comparison ${v\to f}$. Using that $B$ is $\b$-perfect, passing to the perfect subobject as in Lemma~\ref{Lemma-Commutator-Perfect} gives us the needed universal $\b$-central extension $u\colon {U\to A}$ together with the induced comparison morphism ${v\to f}$. This morphism is still a regular epimorphism by the Birkhoff property of $\B$ (see Subsection~\ref{Subsection-Commutators-and-Centralisation}).
\end{proof}

\subsection{An ``absolute'' property of relative universal central extensions}
It is worth remarking here that a universal $\b$-central extension is always central in an absolute sense, namely, with respect to the abelianisation functor $\ab\colon{\A\to \Ab(\A)}$. Here $\Ab(\A)$ is the Birkhoff subcategory of $\A$ consisting of all objects that admit an internal abelian group structure; see, for instance, \cite{Bourn-Gran}.

\begin{proposition}\label{Proposition-Abelian-Result}
Let $\A$ be a semi-abelian category and $\B$ a Birkhoff subcategory of $\A$. If $f\colon {B\to A}$ is a $\b$-central extension with a $\b$-perfect domain~$B$, then $f$ is $\ab$-central. In particular, universal $\b$-central extensions are always $\ab$-central.
\end{proposition}
\begin{proof}
We have $B\cong [B,B]_{\b}$ since $B$ is $\b$-perfect and $[B,B]_{\b}\cong [{B\times_A B},{B\times_A B}]_{\b}$ because $f$ is $\b$-central. Hence the diagonal ${B\to B\times_A B}$, being isomorphic to $\mu_{B}\colon{[B\times_A B,B\times_A B]_{\b}\to B\times_A B}$, is a kernel. By Proposition~3.1 in~\cite{Bourn-Gran}, this implies that $f$ is $\ab$-central. Finally, if $f\colon{B\to A}$ is a universal $\b$-central extension then $B$ is $\b$-perfect.
\end{proof}

\section{Nested Birkhoff subcategories}\label{Section-Nested-Birkhoff-Subcategories}
We now consider the situation where a Birkhoff subcategory $\B$ of a semi-abelian category $\A$ has a further Birkhoff subcategory $\C$ so that they form a chain of nested semi-abelian categories with enough projectives, $\C\subset\B\subset\A$. For instance, $\C$ could be $\Ab(\A)$ as in Theorem~\ref{Theorem-Homology} below. Then there is a commutative triangle of left adjoint functors (all right adjoints are inclusions):
\[
\xymatrix@!0@=4.5em{{\A} \ar[rr]^-{\b} \ar[dr]_-{\cc\b} && \B \ar[dl]^-{\cc}\\
& \C}
\]
Since the objects and morphisms of $\B$ are also objects and morphisms of $\A$, it is natural to compare the notions of $\cc$-centrality, $\cc$-perfect object, homology with respect to $\cc$, etc.\ with that of $\cc\b$-centrality, $\cc\b$-perfect object or the homology with respect to $\cc\b$. We obtain a short exact sequence which relates the two induced types of universal central extension.

\begin{lemma}\label{Lemma-Comparison}
Under the given circumstances:
\begin{enumerate}
\item an object of $\B$ is $\cc$-perfect if and only if it is $\cc\b$-perfect;
\item an extension in $\B$ is $\cc$-central if and only if it is $\cc\b$-central;
\item an extension of $\A$ is $\b$-central as soon as it is $\cc\b$-central.
\end{enumerate}
\end{lemma}
\begin{proof}
If $B$ is an object of $\B$ then $\cc(B)=\cc\b(B)$, which proves the first statement. As for the second statement, an extension $f\colon{B\to A}$ in $\B$ is $\cc$-central if and only if the square in the diagram
\[
\xymatrix{B\times_A B \ar@{-{ >>}}[r]^-{f_{0}} \ar@{-{ >>}}[d]_-{\eta^{\cc}_{B\times_A B}} & B \ar@{-{ >>}}[r]^-{f} \ar@{-{ >>}}[d]^-{\eta^{\cc}_{B}} & A\\
\cc(B\times_A B) \ar@{-{ >>}}[r]_-{\cc(f_{0})} & \cc(B)}
\]
in $\B$ is a pullback. Now the inclusion of $\B$ into $\A$ preserves and reflects all limits and moreover $\cc(f_{0})=\cc\b(f_{0})$, so that $f$ being $\cc$-central is equivalent to $f$ being $\cc\b$-central. The third statement follows from the fact that $\b$ preserves the pullback
\[
\xymatrix{B\times_A B \ar@{-{ >>}}[r]^-{f_{0}} \ar@{-{ >>}}[d]_-{\eta^{\cc\b}_{B\times_A B}} & B \ar@{-{ >>}}[d]^-{\eta^{\cc\b}_{B}}\\
\cc\b(B\times_A B) \ar@{-{ >>}}[r]_-{\cc\b(f_{0})} & \cc\b(B)}
\]
in $\A$ for any $\cc\b$-central extension $f$ in $\A$.
\end{proof}

\begin{lemma}\label{Lemma-Adjunction}
For any object $B$ of $\B$, the adjunction~\eqref{Birkhoff-Adjunction} restricts to an adjunction
\[
\xymatrix@1{{\Centr_{\cc\b}(B)} \ar@<1ex>[r]^-{\b} \ar@{}[r]|-{\perp} & {\Centr_{\cc}(B).} \ar@<1ex>[l]^-{\supset}}
\]
Hence the functor $\b$ preserves universal central extensions:
\[
\b\bigl(u^{\cc\b}_{B}\colon \U(B,\cc\b)\to B\bigr)\,\cong\, \bigl(u^{\cc}_{B}\colon \U(B,\cc)\to B\bigr),
\]
for any $\cc$-perfect object $B$.
\end{lemma}
\begin{proof}
First of all, by Lemma~\ref{Lemma-Comparison}.2, $\Centr_{\cc}(B)\subset (\B\downarrow B)$ is a subcategory of $\Centr_{\cc\b}(B)\subset (\A\downarrow B)$.

Suppose that $g\colon{C\to B}$ is a $\cc\b$-central extension. Applying the functor $\b$, we obtain the extension $\b(g)=g\comp\eta_{C}^{\b}\colon{\b(C)\to B}$, which is $\cc\b$-central as a quotient of~$g$. Being an extension in $\B$, $\b(g)$ is $\cc$-central by Lemma~\ref{Lemma-Comparison}.2.

Finally, as any left adjoint functor, $\b$ preserves initial objects.
\end{proof}

\begin{proposition}\label{Proposition-Comparison}
Suppose that $\C\subset \B\subset \A$ is a chain of inclusions of Birkhoff subcategories of a semi-abelian category $\A$. If $B$ is a $\cc$-perfect object of $\B$ then we have the exact sequence
\[
\xymatrix{0 \ar[r] & [\U(B,\cc\b),\U(B,\cc\b)]_{\b} \ar@{{ |>}->}[r] & \H_{2}(B,\cc\b)\ar@{ >>}[r] & \H_{2}(B,\cc) \ar[r] & 0,}
\]
and $u^{\cc\b}_{B}=u^{\cc}_{B}$ if and only if $[\U(B,\cc\b),\U(B,\cc\b)]_{\b}$ is zero.
\end{proposition}
\begin{proof}
By Lemma~\ref{Lemma-Adjunction} and Theorem~\ref{Theorem-Universal-Central-Extension}, if $B$ is a $\cc$-perfect object of $\B$ then the comparison morphism between the induced universal central extensions gives rise to the $3\times 3$ diagram in Figure~\ref{3x3}.
\begin{figure}
\resizebox{.8\textwidth}{!}{
\xymatrix{& 0 \ar[d] & 0 \ar[d] \\
& [\U(B,\cc\b),\U(B,\cc\b)]_{\b} \ar@{{ |>}->}[d] \ar@{=}[r] & [\U(B,\cc\b),\U(B,\cc\b)]_{\b} \ar@{{ |>}->}[d]\ar@{ >>}[r] & 0 \ar@{{ |>}->}[d] \\
0 \ar[r] & \H_{2}(B,\cc\b) \ar@{ >>}[d] \ar@{{ |>}->}[r] & \U(B,\cc\b) \ar@{ >>}[d]_-{\eta_{\U(B,\cc\b)}^{\b}} \ar@{ >>}[r]^-{u^{\cc\b}_{B}} & B \ar@{=}[d] \ar[r] & 0\\
0 \ar[r] & \H_{2}(B,\cc) \ar@{{ |>}->}[r] \ar[d] & \U(B,\cc) \ar@{ >>}[r]_-{u^{\cc}_{B}} \ar[d] & B \ar[r] \ar[d] & 0\\
& 0 & 0 & 0}}
\caption{The proof of Proposition~\ref{Proposition-Comparison}}\label{3x3}
\end{figure}
\end{proof}

\section{Recognition of universal central extensions}\label{Section-Recognition}
We now prove some recognition results on universal $\b$-central extensions. For these results to hold, we shall need that in $\A$, composition of $\ab$-central extensions is well-behaved, in the sense that any composite
\[
\xymatrix{C \ar@{-{ >>}}[r]^g & B \ar@{-{ >>}}[r]^f & A}
\]
of $\ab$-central extensions $f$ and $g$ with $B$ an $\ab$-perfect object is again $\ab$-central. We shall also ask that $\Ab(\A)$ is contained in $\B$, so that we may suitably reduce the given relative situation to the absolute case. The examples given in~\ref{Subsection-Counterexamples} explain why these conditions are not automatically satisfied. The main result we work towards is Theorem~\ref{Theorem-Homology}, which says that a $\b$-central extension $u\colon U\to A$ is universal if and only if $\H_{1}(U,\b)$ and $\H_{2}(U,\b)$ are zero.

\begin{definition}\label{Condition-U}
Let $\A$ be a semi-abelian category with enough projectives. We say that $\A$ satisfies \defn{\UCEspecial} when the following holds: if $B$ is an $\ab$-perfect object and $f\colon{B\to A}$ and $g\colon{C\to B}$ are $\ab$-central extensions then the extension $f\comp g$ is $\ab$-central.
\end{definition}

\begin{lemma}\label{Lemma-Composition-Central-Extensions-Relative}
Let $\A$ be a semi-abelian category with enough projectives satisfying \UCE\ and~$\B$ a Birkhoff subcategory of $\A$ that contains $\Ab(\A)$. If $u\colon{U\to A}$ is a $\b$-central extension and $v\colon{V\to U}$ is a universal $\b$-central extension then the extension $u\comp v$ is $\b$-central.
\end{lemma}
\begin{proof}
By Proposition~\ref{Proposition-Abelian-Result} both $u$ and $v$ are $\ab$-central. Moreover, as $\Ab(\A)$ is contained in the Birkhoff subcategory~$\B$ of~$\A$, the objects $U$, $V$ and $A$ are $\ab$-perfect. Now by \UCE, also the composite $u\comp v\colon{V\to A}$ is $\ab$-central. Again using that $\B$ is bigger than $\Ab(\A)$ we see that $u\comp v\colon{V\to A}$ is a $\b$-central extension (cf.\ Lemma~\ref{Lemma-Comparison}.3).
\end{proof}

Under the given assumptions, $u\comp v$ is in fact universal, as shown in Proposition~\ref{Proposition-Composition-Universal-Central-Extensions}.

\begin{proposition}\label{Proposition-Universal-vs-Weak-Universal'}
Let $\A$ be a semi-abelian category with enough projectives satisfying \UCE\ and~$\B$ a Birkhoff subcategory of $\A$ that contains $\Ab(\A)$. Then in Proposition~\ref{Proposition-Universal-vs-Weak-Universal}, condition~4 implies condition 1. Hence a~$\b$-central extension $u\colon{U\to A}$ is universal if and only if its domain $U$ is $\b$-perfect and projective with respect to all $\b$-central extensions.
\end{proposition}
\begin{proof}
Suppose that $u\colon{U\to A}$ is a universal $\b$-central extension; we have to prove that every $\b$-central extension of $U$ splits. By Theorem~\ref{Theorem-Universal-Central-Extension}, $U$ admits a universal $\b$-central extension $v\colon{V\to U}$. It suffices to prove that this $v$ is a split epimorphism. By Lemma~\ref{Lemma-Composition-Central-Extensions-Relative}, the composite $u\comp v$ is $\b$-central. The weak $\b$-universality of $u$ now yields a morphism $s\colon{U\to V}$ such that $u\comp v \comp s = u$. But also $u\comp 1_{U}=u$, so that $v\comp s=1_{U}$ by the $\b$-universality of $u$, and the universal $\b$-central extension $v$ splits. The result follows.
\end{proof}

\begin{theorem}\label{Theorem-Homology}
Let $\A$ be a semi-abelian category with enough projectives satisfying \UCE\ and $\B$ a Birkhoff subcategory of $\A$ containing~$\Ab(\A)$. A $\b$-central extension $u\colon U\to A$ is universal if and only if $\H_{1}(U,\b)$ and $\H_{2}(U,\b)$ are zero.
\end{theorem}
\begin{proof}
$\Rightarrow$ If $u\colon U\to A$ is a universal $\b$-central extension then by Proposition~\ref{Proposition-Universal-vs-Weak-Universal'} we have $\H_{1}(U,\b)=\b(U)=0$ and $U$ is projective with respect to all $\b$-central extensions. This implies that $1_{U}\colon{U\to U}$ is a universal $\b$-central extension of~$U$. Theorem~\ref{Theorem-Universal-Central-Extension} now tells us that $\H_{2}(U,\b)=0$.

$\Leftarrow$ The object $U$ is $\b$-perfect because $\b(U) = \H_{1}(U,\b)=0$; since $\H_{2}(U,\b)$ is also zero, the universal $\b$-central extension $u_{U}^{\b}\colon{\U(U,\b)\to U}$ of $U$ induced by Theorem~\ref{Theorem-Universal-Central-Extension} is an isomorphism. Proposition~\ref{Proposition-Universal-vs-Weak-Universal'} now implies that $U\cong \U(U,\b)$ is projective with respect to all $\b$-central extensions. Another application of Proposition~\ref{Proposition-Universal-vs-Weak-Universal'} shows that $u$ is also a universal $\b$-central extension.
\end{proof}

\begin{proposition}\label{Proposition-Composition-Universal-Central-Extensions}
Let $\A$ be a semi-abelian category with enough projectives satisfying \UCE\ and~$\B$ a Birkhoff subcategory of $\A$ that contains $\Ab(\A)$. Let $f\colon{B\to A}$ and $g\colon{C\to B}$ be $\b$-central extensions. Then the composite $f\comp g$ is a universal $\b$-central extension if and only if $g$ is.
\end{proposition}
\begin{proof}
First note that when $g$ is a universal $\b$-central extension then $f\comp g$ is $\b$-central by Lemma~\ref{Lemma-Composition-Central-Extensions-Relative}. The central extensions $f\comp g$ and $g$ have the same domain, and by Proposition~\ref{Proposition-Universal-vs-Weak-Universal'} their universality only depends on a property of this domain.
\end{proof}

Proposition~\ref{Proposition-Universal-vs-Weak-Universal'} has the following partial converse, which shows that in some sense \UCE\ is necessary: if we want that conditions 1--5 in Proposition~\ref{Proposition-Universal-vs-Weak-Universal} are equivalent, independently of the chosen Birkhoff subcategory $\B$ with $\Ab(\A)\subset\B\subset\A$, then the category $\A$ \emph{must} satisfy \UCE. See also Remark~\ref{Remark-Counterexample}.

\begin{proposition}\label{Proposition-Universal-vs-Weak-Universal-Converse}
Let $\A$ be a semi-abelian category with enough projectives and~$\B$ a Birkhoff subcategory of $\A$ that contains $\Ab(\A)$. If in Proposition~\ref{Proposition-Universal-vs-Weak-Universal} all given conditions are equivalent, then the following holds: if $B$ is a $\b$-perfect object and $f\colon{B\to A}$ and $g\colon{C\to B}$ are $\b$-central extensions then the extension $f\comp g$ is $\b$-central. If, in particular, this happens when $\B=\Ab(\A)$, then the category $\A$ satisfies \UCE.
\end{proposition}
\begin{proof}
Let $B$ be a $\b$-perfect object and $f\colon{B\to A}$ and $g\colon{C\to B}$ $\b$-central extensions. Let $u\colon {B\to A}$ be a universal $\b$-central extension. Then by Proposition~\ref{Proposition-Quotient-of-UCE} the uniquely induced comparison morphism $\overline{f}\colon{u\to f}$ is a regular epimorphism. Pulling back $g$ along $\overline{f}$ as in
\[
\xymatrix{C\times_{B}U \ar@{}[rd]|<<<{\pullback} \ar@{-{ >>}}[d]_-{\underline{f}} \ar@<.5ex>@{-{ >>}}[r]^-{\overline{g}} & U \ar@<.5ex>@{.>}[l] \ar@{-{ >>}}[r]^-{u} \ar@{-{ >>}}[d]^-{\overline{f}} & A \ar@{=}[d] \\
C \ar@{-{ >>}}[r]_-{g} & B \ar@{-{ >>}}[r]_-{f} & A}
\]
we obtain a splitting for $\overline{g}$ through Proposition~\ref{Proposition-Universal-vs-Weak-Universal}. Now the composite $u\comp \overline{g}$ is a $\b$-central extension because its pullback $u^{*}(u\comp \overline{g})$ along $u$ is a $\b$-trivial extension, as a composite of two $\b$-trivial extensions (Remark~\ref{Remark-Central-and-Trivial}). Since $\underline{f}$ is a regular epimorphism by regularity of $\A$, the composite $f\comp g$ is a quotient of the $\b$-central extension~$u\comp \overline{g}$, hence is also $\b$-central.
\end{proof}

\UCEcapital\ allows us to obtain the following refinement of Proposition~\ref{Proposition-Comparison}.

\begin{proposition}\label{Proposition-Comparison-Plus}
Suppose that $\Ab(\A)\subset \C\subset \B\subset \A$ is a chain of inclusions of Birkhoff subcategories of a semi-abelian category $\A$ which satisfies \UCE. If $B$ is a $\cc$-perfect object of $\B$ then $u^{\cc\b}_{\U(B,\cc)}=\eta_{\U(B,\cc\b)}^{\b}$ and
\[
[\U(B,\cc\b),\U(B,\cc\b)]_{\b}\cong\H_{2}(\U(B,\cc),\cc\b).
\]
Hence $u^{\cc\b}_{B}=u^{\cc}_{B}$ iff $\eta_{\U(B,\cc\b)}^{\b}$ is an isomorphism iff $\H_{2}(\U(B,\cc),\cc\b)$ is zero.
\end{proposition}
\begin{proof}
In view of Remark~\ref{Remark-Subobject}, $\eta_{\U(B,\cc\b)}^{\b}$ is $\cc\b$-central as a subobject of $u^{\cc\b}_{B}$; hence Proposition~\ref{Proposition-Composition-Universal-Central-Extensions} implies that it is a universal $\cc\b$-central extension of~$\U(B,\cc)$. The result now follows from Theorem~\ref{Theorem-Universal-Central-Extension}.
\end{proof}

\section{Examples}\label{Section-Examples}
In this section we illustrate the theory with some classical and contemporary examples. All categories we shall be considering here are (equivalent to) varieties of $\Omega$-groups, and as such are semi-abelian with enough projectives. As an illustration of Section~\ref{Section-Nested-Birkhoff-Subcategories}, we shall consider the categories $\Gp$ of groups and $\Ab$ of abelian groups; $\Leibniz_{\K}$, $\Lie_{\K}$ and $\Vect_{\K}$ of Leibniz algebras, Lie algebras and vector spaces over a field~$\K$; and the categories $\PXMod$, $\XMod$ and $\Ab\XMod$ of precrossed modules, crossed modules and abelian crossed modules. Then in~\ref{Subsection-Counterexamples} we give two counterexamples which further justify the conditions that appear in Section~\ref{Section-Recognition}.

\subsection{Groups and abelian groups}
The case of groups and abelian groups is well-known and entirely classical, but we think it is worth repeating. The left adjoint $\ab\colon{\Gp\to \Ab}$ to the inclusion of $\Ab$ in $\Gp$ is called the \defn{abelianisation functor}; it sends a group $G$ to its abelianisation $G/[G,G]$. A surjective group homomorphism $f\colon{B\to A}$ is a central extension if and only if the commutator $[\Ker(f),B]_{\ab}=[\Ker(f),B]$ is trivial; given a group $G$ and a normal subgroup $N$ of~$G$, their commutator $[N,G]$ is the normal subgroup of $G$ generated by the elements $[n,g]=ngn^{-1}g^{-1}$ for all $n\in N$ and $g\in G$. Equivalently $f$ is central if and only if the kernel $\Ker(f)$ is contained in the centre
\[
Z(B)=\{z\in B\mid\text{$[z,b]=1$ for all $b\in B$}\}
\]
of $B$. A group $G$ is perfect if and only if it is equal to its commutator subgroup~$[G, G]$. Computing the second integral homology group $\H_{2}(G,\Z)=\H_{2}(G,\ab)$ of a perfect group $G$ is particularly simple: take the universal central extension
$
u_{G}^{\ab}\colon{\U(G,\ab)\to G};
$
its kernel $\Ker(u^{\ab}_{G})$ is $\H_{2}(G,\Z)$. It is well known that $\Gp$ satisfies \UCE.

\subsection{Leibniz algebras, Lie algebras and vector spaces}
Recall~\cite{Loday-Leibniz} that a \defn{Leibniz algebra} $\Lieg$ is a vector space over a field $\K$ equipped with a bilinear operation $[\cdot , \cdot]\colon \Lieg \times \Lieg \to \Lieg$ that satisfies the \defn{Leibniz identity}
\[
[x,[y,z]]=[[x,y],z]-[[x,z],y]
\]
for all $x$, $y$, $z\in \Lieg$. If $[x,x]=0$ for all $x\in\Lieg$ then the bracket is skew-symmetric and the Leibniz identity is the Jacobi identity, so $\Lieg$ is a Lie algebra.

Here there are three inclusions of Birkhoff subcategories, of which the left adjoints form the commutative triangle
\[
\xymatrix@!0@=4.5em{{\Leibniz_{\K}} \ar[rr]^-{\lie} \ar[dr]_-{\vect\circ \lie} && \Lie_{\K} \ar[dl]^-{\vect}\\
& \Vect_{\K}}
\]
The left adjoint $\lie\colon{\Leibniz_{\K}\to \Lie_{\K}}$ (which is usually called the \defn{Liesation functor}) takes a Leibniz algebra $\Lieg$ and maps it to the quotient $\Lieg/\Lieg^{\ann}$, where $\Lieg^{\ann}$ is the two-sided ideal (= normal subalgebra) of $\Lieg$ generated by all elements $[x,x]$ for $x\in \Lieg$. The category $\Vect_{\K}$ may be considered as a subvariety of $\Lie_{\K}$ by equipping a vector space with the trivial Lie bracket; the left adjoint $\vect\colon{\Lie_{\K}\to \Vect_{\K}}$ to the inclusion $\Vect_{\K}\subset \Lie_{\K}$ takes a Lie algebra $\Lieg$ and maps it to the quotient $\Lieg/[\Lieg,\Lieg]$, where $[\Lieg,\Lieg]$ is generated by the elements $[x,y]\in \Lieg$ for all $x$, $y\in\Lieg$.

The notion of central extension obtained in the case of $\Lie_{\K}$ vs.\ $\Vect_{\K}$ is the ordinary notion of central extension of Lie algebras, where the kernel $\Ker(f)$ of~$f\colon{\Lieb\to \Liea}$ should be included in the centre
\[
Z(\Lieb)=\{z\in \Lieb\mid\text{$[z,b]=0$ for all $b\in \Lieb$}\}
\]
of $\Lieb$. Examples of universal $\vect$-central extensions of Leibniz algebras over a field~$\K$ may be found in~\cite{Casas-Ladra}; in this case, the notion of perfect object is the classical one. The categories $\Lie_{\K}$ and $\Leibniz_{\K}$ are well known to satisfy \UCE\: see, for instance, \cite{Gn} for the Leibniz case.

On the other hand, a Leibniz algebra $\Lieg$ is perfect with respect to $\Lie_{\K}$ if and only if $\Lieg=\Lieg^{\ann}$. Moreover, given a Leibniz algebra $\Lieg$, we may consider the two-sided ideal generated by
\[
\{ z\in \Lieg\mid\text{$[g,z]=-[z,g]$ for all $g\in \Lieg$}\};
\]
we call it the \defn{$\Lie_{\K}$-centre} of $\Lieg$ and denote it by $Z_{\Lie}(\Lieg)$. When $\K$ is a field of characteristic different from $2$ then this relative centre allows us to characterise the $\lie$-central extensions of Leibniz algebras over $\K$.

\begin{proposition}\label{Proposition-Leibniz-Lie}
Suppose $\K$ is a field of characteristic different from $2$. For an extension $f\colon{\Lieb\to \Liea}$ of Leibniz algebras over $\K$, the following three conditions are equivalent:
\begin{enumerate}
\item $f\colon{\Lieb\to \Liea}$ is central with respect to $\Lie_{\K}$;
\item $(\Lieb\times_\Liea \Lieb)^{\ann}\cong \Lieb^{\ann}$;
\item $\Ker(f)\leq Z_{\Lie}(\Lieb)$.
\end{enumerate}
\end{proposition}
\begin{proof}
Condition 1 is equivalent to 2 by definition. Now suppose that 2 holds and consider $k\in \Ker(f)$ and $b\in \Lieb$. Then both $[(k,0),(k,0)]=([k,k],[0,0])$ and $[(b-k,b),(b-k,b)] = ([b-k,b-k],[b,b])$ are in $(\Lieb\times_\Liea \Lieb)^{\ann}$, which implies that $[k,k]=[0,0]=0$ and $[b-k,b-k]=[b,b]$. Thus we see that $[b,k] + [k,b] = 0$, which implies that condition 3 holds.

Conversely, consider $[(b-k,b),(b-k,b)]$ in $(\Lieb\times_\Liea \Lieb)^{\ann}\meet \Ker(f_{1})$, where $f_{1}$ denotes the second projection of the kernel pair of $f$. Then $k$ is an element of the kernel of~$f$ and
\[
0=f_{1}[(b-k,b),(b-k,b)]=f_{1}([b-k,b-k],[b,b])=[b,b].
\]
Now 3 implies that $[k,k]+[k,k]=0$, so that $[k,k]=0$, since $\chara(\K)\neq 2$. Furthermore, $[b,k]+[k,b]=0$, which implies that
\[
[b-k,b-k]=[b,b]-[b,k]-[k,b]+[k,k]=0.
\]
Hence $(\Lieb\times_\Liea \Lieb)^{\ann}\meet \Ker(f_{1})$ is zero, so $(\Lieb\times_\Liea \Lieb)^{\ann}\cong \Lieb^{\ann}$ and condition 2 holds.
\end{proof}

Given a Leibniz algebra $\Lieg$, the homology vector space $\H_{2}(\Lieg,\vect\comp \lie)$ is the Leibniz homology developed in~\cite{LP}; see also \cite{CP, Pirashvili-Leibniz}. As far as we know, the homology Lie algebra $\H_{2}(\Lieg,\lie)$ has not been studied before, but certainly the theories referred to in Remark~\ref{Remark-Homology} apply to it. If $\Lieg$ is a Lie algebra then the vector space $\H_{2}(\Lieg,\vect)$ is the classical Chevalley--Eilenberg homology. If we interpret Proposition~\ref{Proposition-Comparison-Plus} in the context of Lie algebras then we recover Corollary~2.7 from~\cite{Gn}, but in the special case where $\K$ is a field:

\begin{proposition}\label{Proposition-Gn}
If $\Lieg$ is a perfect Lie algebra then $\H_{2}(\U(\Lieg,\vect),\vect\comp\lie)$ is the kernel of ${\H_{2}(\Lieg,\vect\comp\lie)\to \H_{2}(\Lieg,\vect)}$. Moreover, 
\[
[\U(\Lieg,\vect\comp\lie),\U(\Lieg,\vect\comp\lie)]_{\lie}=\H_{2}(\U(\Lieg,\vect),\vect\comp\lie).
\]
\end{proposition}
\begin{proof}
By Proposition~\ref{Proposition-Comparison-Plus} it suffices to recall that the category $\Lie_{\K}$ satisfies \UCE.
\end{proof}

\subsection{Precrossed modules, crossed modules and abelian crossed modules}
Recall that a \defn{precrossed module} $(T,G,\del)$ is a group homomorphism $\del\colon{T\to G}$ together with an action of $G$ on $T$, denoted by ${}^{g}t$ for $g\in G$ and $t\in T$, satisfying $\del({}^{g}t)=g\del(t)g^{-1}$ for all $g\in G$ and $t\in T$. If in addition it verifies the \defn{Peiffer identity} ${}^{\del(t)}t'=tt't^{-1}$ for all $t$, $t'\in T$ then we say that $(T,G,\del)$ is a \defn{crossed module}. A morphism of (pre)crossed modules $(f_{1},f_{0})\colon{(T,G,\del)\to (T',G',\del')}$ consists of group homomorphisms $f_{1}\colon{T\to T'}$ and $f_{0}\colon{G\to G'}$ such that $\del'\comp f_{1}=f_{0}\comp \del$ and the action is preserved. The categories $\PXMod$ and $\XMod$ are equivalent to varieties of $\Omega$-groups; see, for instance, \cite{Janelidze-Marki-Tholen}, \cite{LR} or \cite{Loday}. The category $\Ab\XMod$ consists of \defn{abelian} crossed modules: $(T,G,\del)$ such that $T$ and $G$ are abelian groups and the action of $G$ on $T$ is trivial.

As for algebras, we obtain a commutative triangle of left adjoint functors.
\[
\xymatrix@!0@=4.5em{{\PXMod} \ar[rr]^-{\xmod} \ar[dr]_-{\ab\circ \xmod} && \XMod \ar[dl]^-{\ab}\\
& \Ab\XMod}
\]
Given two normal precrossed submodules $(M,H,\del)$ and $(N,K,\del)$ of a precrossed module $(T,G,\del)$, the \defn{Peiffer commutator} $\langle M,N\rangle$ is the normal subgroup of~$T$ generated by the Peiffer elements
\[
\langle m,n\rangle=mnm^{-1}({}^{\del(m)}n)^{-1}
\qquad\text{and}\qquad
\langle n,m\rangle=nmn^{-1}({}^{\del(n)}m)^{-1}
\]
for $m\in M$, $n\in N$~\cite{El}. We denote by $\langle (M,H,\del),(N,K,\del)\rangle$ the precrossed module $(\langle M,N\rangle, 0,0)$; it may be considered as a normal precrossed submodule of $(T,G,\del)$. The precrossed module $\langle(T,G,\del),(T,G,\del)\rangle =(\langle T,T\rangle, 0, 0)$ is the smallest one that makes the quotient $(T,G,\del)/\langle(T,G,\del),(T,G,\del)\rangle$ a crossed module. This defines a functor $\xmod\colon{\PXMod\to \XMod}$, left adjoint to the inclusion of $\XMod$ in $\PXMod$.

A precrossed module $(T,G,\del)$ is $\xmod$-perfect if and only if
\[
\langle (T,G,\del),(T,G,\del)\rangle=(T,G,\del).
\]
In particular, then $G=0$; hence $\langle T,T\rangle=[T,T]$, so that $(T,G,\del)$ is perfect with respect to $\XMod$ exactly when $T$ is perfect with respect to $\Ab$ and $G$ is trivial.

The results of~\cite[Section~9.5]{EGVdL} imply that an extension of precrossed modules $f\colon{B\to A}$ is central with respect to $\XMod$ if and only if $\langle \Ker(f),B\rangle=1$; the following characterisation may also be shown analogously to Proposition~\ref{Proposition-Leibniz-Lie}. Given a precrossed module $(T,G,\del)$, its \defn{$\XMod$-centre} $Z_{\XMod}(T,G,\del)$ is the normal precrossed submodule $(Z_{\XMod}(T),G,\del)$ of $(T,G,\del)$ where
\[
Z_{\XMod}(T)= \{t\in T\mid\text{$\langle t,t'\rangle=1=\langle t',t\rangle$ for all $t'\in T$}\}.
\]

\begin{proposition}
For an extension $(f_{1},f_{0})\colon{(T,G,\del)\to (T',G',\del')}$ of pre\-cros\-sed modules, the following conditions are equivalent:
\begin{enumerate}
\item $(f_{1},f_{0})$ is central with respect to $\XMod$;
\item $\langle (T\times_{T'}T,G\times_{G'}G, \del\times_{\del'}\del),(T\times_{T'}T,G\times_{G'}G, \del\times_{\del'}\del)\rangle\cong \langle (T,G,\del),(T,G,\del)\rangle$;
\item $\langle T\times_{T'}T,T\times_{T'}T\rangle\cong \langle T,T\rangle$;
\item $\Ker(f_{1})\leq Z_{\XMod}(T)$;
\item $\Ker(f_{1},f_{0})\leq Z_{\XMod}(T,G,\del)$.\hfill\qed
\end{enumerate}
\end{proposition}

We now focus on the further adjunction to $\Ab\XMod$. Given a precrossed module $(T,G,\del)$, the commutator $[G,T]$ is the normal subgroup of $T$ generated by the elements ${}^{g}tt^{-1}$ for $g\in G$ and $t\in T$. The left adjoint $\ab\colon{\XMod \to \Ab\XMod}$ takes a crossed module $(T,G,\del)$ and maps it to $\ab(T,G,\del)=(T/[G,T],G/[G,G],\overline{\del})$, where $\overline{\del}$ is the induced group homomorphism. The functor
\[
\ab\comp \xmod\colon{\PXMod \to \Ab\XMod}
\]
maps a precrossed module $(T,G,\del)$ to $(T/[T,T][G,T],G/[G,G],\overline{\del})$. 

As shown in~\cite{Bourn-Gran}, an extension of crossed modules is central with respect to $\Ab\XMod$ exactly when it is central in the sense of~\cite{GL}. An extension of \emph{pre}crossed modules is central with respect to $\Ab\XMod$ if and only if it is central in the sense of~\cite{ACL, AL}. In this case, the notions of perfect object obtained are classical. The article~\cite{ACL} gives several non-trivial examples of universal central extensions of (pre)crossed modules, relative to $\Ab\XMod$. Lemma~13 in~\cite{ALG2} implies that, for abelianisation of precrossed modules, all the conditions in Proposition~\ref{Proposition-Universal-vs-Weak-Universal} are equivalent. Hence by Proposition~\ref{Proposition-Universal-vs-Weak-Universal-Converse}, \UCE\ holds in $\PXMod$. Lemma~\ref{Lemma-Comparison} now tells us that also $\XMod$ satisfies \UCE.

The homology crossed module $\H_{2}((T,G,\del),\ab\comp \xmod)$ was studied in~\cite{ALG}, while $\H_{2}((T,G,\del),\ab)$ was considered in~\cite{Carrasco-Homology}. For a precrossed module $(T,G,\del)$, the relative $\H_{2}((T,G,\del),\xmod)$ was characterised in~\cite{EGVdL}. If we interpret Proposition~\ref{Proposition-Comparison} in this situation then we regain~\cite[Theorem~5]{ACL}.

\subsection{Two counterexamples}\label{Subsection-Counterexamples}
Our first counterexample is borrowed from~\cite{CIP}. It shows that a category---here the category $\NAAlg_\K$ of non-associative algebras over a field $\K$, which is a variety of $\Omega$-groups---can be semi-abelian without having to satisfy \UCE. This means that $\NAAlg_\K$ does not quite match the picture sketched in Section~\ref{Section-Recognition}.

We must also emphasise that \UCE\ by itself is not yet strong enough to yield results such as Proposition~\ref{Proposition-Universal-vs-Weak-Universal'} or Theorem~\ref{Theorem-Homology}, unless $\Ab(\A)$ is contained in $\B$. Example~\ref{Counter2}, which explains this, was offered to us by George Peschke. It describes a universal~$\b$-central extension $u\colon{U\to A}$ such that $\H_{2}(U,\b)$ is not trivial---and indeed one of the assumptions of Theorem~\ref{Theorem-Homology} is violated: the Birkhoff subcategory~$\B$ of the \mbox{(semi-)}abelian category $\A$ which we shall consider is strictly smaller than $\Ab(\A)$.

\begin{example}\label{Counter1}
A \defn{non-associative algebra} $\Lieg$ is a vector space over a field $\K$ equipped with a bilinear operation $[\cdot , \cdot]\colon \Lieg \times \Lieg \to \Lieg$. Unlike for Leibniz or Lie algebras (or for associative ones), the bracket need not satisfy any additional conditions. We write $\NAAlg_\K$ for the category of non-associative algebras over $\K$ and remark that it coincides with the category of $\hom$-Leibniz algebras of which the twisting map is trivial ($\alpha=0$ in the notations of~\cite{MS, CIP}) and with the category of magmas in the monoidal category $(\Vect_{\K}, \tensor, \K)$. Note that $\Leibniz_{\K}$, and hence also $\Lie_\K$ and $\Vect_\K$, are subvarieties of $\NAAlg_\K$. Furthermore, an algebra is abelian if and only if it has a trivial bracket, so precisely when it ``belongs to'' $\Vect_\K$. We write $\vect\colon{\NAAlg_{\K}\to \Vect_{\K}}$ for the abelianisation functor. It is easily seen that an extension $f\colon{\Lieb\to \Liea}$ in $\NAAlg_\K$ is $\vect$-central if and only if its kernel is contained in the centre of $\Lieb$,
$
Z(\Lieb)=\{z\in \Lieb\mid\text{$[z,b]=0=[b,z]$ for all $b\in \Lieb$}\}.
$
 
In~\cite{CIP} the following situation is considered. 
\[
\xymatrix{\Liec \ar@{-{ >>}}[rr]^-{g} \ar@{-{ >>}}[rd]_-{f\circ g} && \Lieb \ar@{-{ >>}}[ld]^-{f}\\
& \Liea}
\]
As vector spaces, $\Liea$, $\Lieb$ and $\Liec$ are $2$-, $3$- and $4$-dimensional with respective bases $\{a_{1},a_{2}\}$, $\{b_{1},b_{2},b_{3}\}$ and $\{c_{1},c_{2},c_{3},c_{4}\}$. Their brackets are generated by
\begin{gather*}
[a_{2},a_{1}]=a_{2},\quad [a_{2},a_{2}]=a_{1}\\
[b_{2},b_{2}]=b_{1},\quad [b_{3},b_{2}]=b_{3}, \quad [b_{3},b_{3}]=b_{2}\\
[c_{3},c_{2}]=c_{1}, \quad [c_{3},c_{3}]=c_{2}, \quad [c_{4},c_{3}]=c_{4} \quad [c_{4},c_{4}]=c_{3}
\end{gather*}
and zero elsewhere. The algebras $\Lieb$ and $\Liea$ are $\vect$-perfect because $[\Lieb,\Lieb]=\Lieb$. The morphism of non-associative algebras $f$ sends $(b_{1},b_{2},b_{3})$ to $(0,a_{1},a_{2})$ and $g$ sends $(c_{1},c_{2},c_{3},c_{4})$ to $(0,b_{1},b_{2},b_{3})$. The kernel of $f$ is generated by $b_{1}$ and thus equal to the centre $Z(\Lieb)$ of $\Lieb$. Hence $f$ is $\vect$-central. Likewise, $\Ker(g)$ is generated by $c_{1}$ and thus equal to $Z(\Liec)$, so that $g$ is $\vect$-central. On the other hand, the kernel of~$f\comp g$ contains $c_{2}$, so that $f\comp g$ cannot be $\vect$-central.
\end{example}

\begin{remark}\label{Remark-Counterexample}
Combining Example~\ref{Counter1} with Proposition~\ref{Proposition-Universal-vs-Weak-Universal-Converse} we obtain a contradiction with the statement of Proposition~6.3 in~\cite{Cheng-Su}. It entails the equivalence of all conditions in Proposition~\ref{Proposition-Universal-vs-Weak-Universal}, which would mean that the category of $\hom$-Leibniz algebras satisfies \UCE. We know, however, that already its subvariety $\NAAlg$ does not, which through Lemma~\ref{Lemma-Comparison} leads to a clash. It appears that the proof given in~\cite{Cheng-Su} does not explain the second half of the ``necessary condition''.
\end{remark}

\begin{example}\label{Counter2}
Let $C$ be the infinite cyclic group (with its generator written $c\in C$) and $R=\Z[C]$ the integral group-ring over $C$. We take $\A$ to be the (abelian) category ${}_{R}\Mod$ of modules over~$R$, so that $\Ab(\A)=\A$ and \UCE\ holds. We consider its full subcategory $\B$ of all $R$-modules with a trivial $C$-action; it is clearly Birkhoff in $\A$, and its reflector is determined by tensoring with the trivial $R$-module $\Z$, so that ${\b(M)=\Z\tensor_{R}M}$ for any $R$-module $M$.

Now consider a prime number $p\neq 2$ and let $M$ be the $R$-module $\bigvee_{k\geq 1}M_{k}$, where~$M_{k}$ for $k\geq 1$ is the abelian group $\Z_{p^{k}}=\Z/_{p^{k}\Z}$ equipped with the $C$-action
\[
c\cdot m=(1-p)\cdot m.
\]
Note that a natural inclusion of $R$-modules $M_{k}\to M_{k+1}$ is given by
\[
(l+p^{k}\Z)\mapsto (p\cdot l+p^{k+1}\Z).
\]
Then it may be checked that $\H_{2}(M,\b)=\H_{2}(C,M)\cong \Z_{p}\neq 0$, while $M$ is $\b$-perfect, and
\[
u\colon M\to M\colon m\mapsto p\cdot m
\]
is a universal $\b$-central extension.
\end{example}

\section{Final remarks}\label{Section-Final-Remarks}
\subsection{When does \UCE\ hold?}
The main question left unanswered in this paper is, under which precise conditions on the category $\A$ the results of Section~\ref{Section-Recognition} are valid. That is to say, are there elementary characterisations for \UCE, or at least convenient sufficient conditions? There is a whole list of categorical conditions which may be relevant to \UCE. One good candidate seems to be the context of semi-abelian \emph{action accessible} categories where the concept of \emph{centraliser} is well-behaved~\cite{BJ07}, but there are also \emph{Moore} categories~\cite{Rodelo:Moore}, \emph{action representable} categories~\cite{BJK, BJK2} and the \emph{Smith is Huq} condition~\cite{MFVdL}, to name a few others. Further note that \UCE\ already makes sense in the context of Section~\ref{Section-Basic-Definitions}: pointed Barr exact Goursat categories.

\subsection{Stem extensions}
The last chapter of Julia Goedecke's thesis~\cite{JuliaThesis} contains results on \emph{stem extensions} and \emph{stem covers} in the absolute context of a semi-abelian category $\A$ with the chosen Birkhoff subcategory $\Ab(\A)$. 
In the process, she obtains a categorical version of results from~\cite{VC} on central extensions of crossed modules. At least part of the theory should also be valid in a relative context.

\subsection{Higher-order extensions}
The article~\cite{Peschke-UE} sketches how universality of high\-er-or\-der extensions may be introduced and used in some concrete situations in algebra and in topology. In the semi-abelian case its results imply that universal \emph{higher} central extensions---in the sense of~\cite{EGVdL}, and relative to abelianisation---exist for objects which are perfect in the appropriate higher-order sense. More explicitly, assuming that enough projectives exist and the \emph{Smith is Huq} condition holds, an object~$A$ will admit a universal $n$-fold central extension precisely when all homology objects $\H_{i}(A,\ab)$ for $1\leq i\leq n$ vanish. Here the word ``universal'' needs to be treated with some care, because higher central extensions may be equivalent (belong to the same cohomology class) without being isomorphic~\cite{RVdL2}.

\section*{Acknowledgements}

We would like to thank 
 Tomas Everaert, Julia Goedecke and George Peschke for some invaluable suggestions. Thanks also to the University of Coimbra, the University of Vigo and the Banff International Research Station for their kind hospitality.

%

\providecommand{\noopsort}[1]{}
\providecommand{\bysame}{\leavevmode\hbox to3em{\hrulefill}\thinspace}
\providecommand{\MR}{\relax\ifhmode\unskip\space\fi MR }
\providecommand{\MRhref}[2]{%
  \href{http://www.ams.org/mathscinet-getitem?mr=#1}{#2}
}
\providecommand{\href}[2]{#2}

\end{document}